\documentclass[12pt]{article}
\usepackage[a4paper,left=2cm,right=2cm,top=2cm,bottom=2cm]{geometry}
\usepackage[utf8]{inputenc}
\usepackage[leqno]{amsmath}
\usepackage{amssymb}
\usepackage[T1]{fontenc}
\usepackage{graphicx}
\usepackage{subcaption}
\usepackage{chngcntr}
\counterwithin{figure}{section}
\counterwithin{equation}{section}
\usepackage{xcolor}
\usepackage{hyperref}
\usepackage[all]{hypcap} 
\hypersetup{
  colorlinks   = true, 
  urlcolor     = blue, 
  linkcolor    = blue, 
  citecolor   = blue 
}
\usepackage{epstopdf, epsfig}
\usepackage[round]{natbib}

\def\R{\mathbb{R}}

\def\v{\boldsymbol{v}}
\def\w{\boldsymbol{w}}
\def\k{\boldsymbol{k}}

\def\a{\boldsymbol{a}}

\def\x{\boldsymbol{x}}
\def\3{\boldsymbol{e_3}}
\def\c{\boldsymbol{c}}
\def\bzeta{\boldsymbol{\zeta}}
\def\bpsi{\boldsymbol{\psi}}
\def\F{\boldsymbol{F}}
\def\tb{{\tilde b}}
\def\kb{{\bar\k}}
\def\kbp{{\bar\k}^{\perp}}
\def\p{\partial}

\newtheorem{remark}{Remark}

\title{Explicit superposed and forced plane wave generalized Beltrami flows}
\author{Artur Prugger \& Jens D. M. Rademacher}
\date{\today}

\begin{document}

\maketitle

We revisit and present new linear spaces of explicit solutions to incompressible Euler and Navier-Stokes equations on $\R^n$, as well as the rotating 
Boussinesq equations on $\R^3$. We cast these solutions are superpositions of certain linear plane waves of arbitrary amplitudes that also solve 
the nonlinear equations by constraints on wave vectors and flow directions. For $n\leq 3$ these are explicit examples for generalized 
Beltrami flows. We show that forcing terms of corresponding plane wave type yield explicit solutions by linear variation of constants. We 
work in Eulerian coordinates and distinguish the two situations of vanishing and of gradient nonlinear terms, where the 
nonlinear terms modify the pressure. The methods, that we introduce to find explicit solutions in some nonlinear fluid models, can also be used in other equations 
with material derivative. Our approach offers another view on known explicit solutions of different fluid models from a plane wave perspective,
and provides transparent nonlinear interactions between different flow components. \newline

\section{Introduction}
Explicit solutions form a cornerstone in the concrete analysis of nonlinear models and continue to be of relevance in fluid models, e.g., \cite{Achatz06,Drazin06,Toorn2019,Dyck19}.
Such solutions provide insight into the mathematical structure of models, can be practical test cases for numerical schemes, and 
organising centers for relevant dynamics. Their vicinity can be analysed by perturbation methods, such as weakly nonlinear analysis, where parameters of a known family of solutions are modulated on selected spatio-temporal and amplitude 
scales. For the resulting reduced modulation equations in spatially extended systems, often the interesting aspects of dynamics are built from wave-like phenomena, so one is interested 
in the existence, stability and interaction of these.

In this paper we revisit some known and present -- to the best of our knowledge -- new families of linear spaces of explicit solutions to Eulerian fluid models 
on the whole space $\R^n$.
Coming from a wave perspective, and for better comparison of solutions as well as investigations of nonlinear interactions between flow components, 
we seek superpositions of travelling or plane waves. These should admit a free scaling and  
simultaneously solve the nonlinear and the linear equations that arise upon dropping the nonlinear term. 
Such solutions do not occur in generic nonlinear evolution equations, but it is well known that these are possible for the nonlinear terms of the material derivative in fluid models.
Simple explicit monochromatic wave solutions of this kind, i.e., consisting of a single Fourier mode, are presented, e.g., in 
\cite{Meshalkin61, Mied76, Drazin77}. More specifically, for $n\leq 3$ the solutions we present can be viewed as superposed explicit 
generalized Beltrami flows \citep{Wang89,Wang90,Drazin06}, which is a class of solutions satisfying the condition
\[
\nabla\times\bigl(\v\times(\nabla\times \v)\bigr)=0\,.
\]
Related to this are the classes of Beltrami flows, satisfying $\v\times(\nabla\times \v)=0$ \citep{Majda02,Wang89,Drazin06}, and 
the extended Beltrami flows \citep{Dyck19}, to which our solutions do not belong.

The basic idea relies on orthogonality of wave vectors and flow directions, for which the nonlinear term $(\v \cdot \nabla) \w$ of the material transport
vanishes, similar to transverse waves. As an illustration, consider the 2D barotropic quasi-geostrophic equation for the horizontal stream function $\psi$ given by
\begin{align*}
\partial_t \Delta\psi + \nabla^{\perp}\psi\cdot\nabla \Delta\psi = 0 \,, \quad \x\in\R^2 \,,
\end{align*}
with $\nabla^{\perp}=(-\partial_y,\partial_x)$ and 
 corresponding linear equation $\partial_t \Delta\psi=0$ (e.g., equation (1.51) in \cite{M2}). Any stationary single mode $\alpha \sin(\k\cdot\x)$ 
with arbitrary wave vector $\k\in\R^2$ and amplitude $\alpha$ solves the linear and also the nonlinear equation due to the 
orthogonality of $\nabla^\perp \psi$ and $\nabla \Delta\psi$. 
However, this is no longer the case for general superpositions of such monochromatic wave modes with different $\k$.

In this paper, we explore this idea further and revisit as well as newly identify a number of admissible superpositions in Eulerian coordinates for the basic 
rotating and non-rotating incompressible Euler and Navier Stokes equations on $\R^n$, as well as the Boussinesq equations on 
$\R^3$, 
all expressed in terms of velocities. Under suitable forcing these solutions yield steady states; an example are the solutions in 
\cite{Meshalkin61} and the Kolmogorov flows discussed in \cite{BalmforthYoung2002,BalmforthYoung2005}, though in this paper we 
will not further discuss stability properties of such states. Other results about superimposed flow can be found in \cite{Kambe86, 
Hui87, Majda02,  Drazin06}, and we remark on the related literature further in the Discussion \S\ref{s:discuss}. 
In contrast to most of these cases, and in line with, e.g., \cite{Walsh}, we also consider nonzero nonlinear term $(\v \cdot \nabla) \w$,
which however is a gradient, such that it contributes to pressure only. 

\medskip
In this paper we want to be explicit in the velocities and pressure, and therefore do not consider vorticity or Lagrangian coordinates. 
We investigate explicit solutions with the plane wave or Fourier mode approach, in order to transparently compare the nonlinear interactions between the flow components. This also simplifies explicitly determining the pressure, which is compensating the nonlinear term and, in the rotating case, also the Coriolis term.
The new contributions in this paper are on the one hand the clarification of the possible superpositions for explicit flows, including different scales as well as Coriolis term. On the other hand the treatment of suitable forcing as linear dynamics, and the (somewhat abstract) consideration of arbitrary spatial dimensions. 

The solutions considered here that come with a gradient nonlinear term are built from planar flows,
where the wave vector directions within each plane are arbitrary, but necessarily of the same wave length, i.e., on the same scale \citep{Hui87, Majda02,Walsh,Chai20}. 
Moreover, we discuss dimensions of parameterised solution spaces and explicitly specify the pressure in each case. 
Regarding related solutions, we present explicit superposed monochromatic gravity waves and Kolmogorow flow, which we have not found elsewhere.

The paper is organised as follows: We discuss the basic approach via the rotating Boussinesq equations in \S\ref{s:bouss}, and turn 
to non-rotating models in \S\ref{s:nonrot} for more general forms of solutions in $\R^n$. In \S\ref{s:forcing} we consider adapted forcing, and end with a discussion in 
\S\ref{s:discuss}.

\section{Basic approach via the rotating Boussinesq equations}\label{s:bouss}
We start with presenting classes of explicit nonlinear plane wave solutions to the rotating Boussinesq equations. These sets of solutions 
are more restricted than in simpler fluid models and larger sets of solutions will be derived based on the same ideas for the 
non-rotating 
case as well as Navier-Stokes and Euler equations.\newline
The unforced viscous rotating Boussinesq equations in the f-plane approximation on $\R^3$ read
\begin{subequations}\label{eq 1}
\begin{align}
\frac{\partial\boldsymbol{v}}{\partial t}+(\boldsymbol{v}\cdot\nabla)\boldsymbol{v}+f\boldsymbol{e_3}\times\boldsymbol{v} + \nabla p - 
\boldsymbol{e_3}b \ &= \ \nu\Delta\boldsymbol{v}\label{eq 1a}\\[2mm]
\nabla\cdot\boldsymbol{v} \ &= \ 0 \label{eq 1b}\\[2mm]
\frac{\partial b}{\partial t}+(\boldsymbol{v}\cdot\nabla)b-\frac{d\bar{\rho}}{dz}v_3 \ &= \ \mu\Delta b\,,\label{eq 1c}
\end{align}
\end{subequations}
with velocity field $\boldsymbol{v}(t,\boldsymbol{x})\in\mathbb{R}^3$ for $\boldsymbol{x}\in\mathbb{R}^3$ and $t\geq 0$, pressure and 
buoyancy $p(t,\boldsymbol{x}), b(t,\boldsymbol{x})\in\mathbb{R}$, the Coriolis parameter $f\in\mathbb{R}\backslash\{0\}$, the vertical 
unit vector $\boldsymbol{e_3}$, viscosity parameter $\nu \geq 0$, thermal diffusivity $\mu\geq 0$, and reference density field 
$\bar{\rho}(z)$, 
where as usual we assume linear dependence of the vertical space direction $z$ for linear stratification. More specifically, the buoyancy 
satisfies $b(t,\boldsymbol{x})=-g(\rho(t,\boldsymbol{x})-\bar{\rho}(z))/\rho_0\in\mathbb{R}$, with fluid density 
$\rho(t,\boldsymbol{x})\in\mathbb{R}$, 
characteristic density $\rho_0$ and gravitational acceleration $g$. See, e.g., \cite{Achatz06, GohWayne}. 
We shall focus on $\nu,\mu>0$, but remark on the inviscid case $\nu=\mu=0$.

\subsection{Vanishing nonlinearity}\label{s:vanish1}
It is well known that the nonlinear term $\v\cdot\nabla\v$ in \eqref{eq 1a} vanishes for pure plane waves
\begin{align*}
\boldsymbol{v}(t,\boldsymbol{x}) = \psi(t,\boldsymbol{k}\cdot\boldsymbol{x})\boldsymbol{a}\,,
\end{align*}
with arbitrary scalar valued function $\psi\in C^{2}(\mathbb{R}_{\geq 0}\times\mathbb{R},\R)$, and orthogonal wave vector and constant flow direction 
\begin{align*}
\k, \a \in\R^3 \,,\quad \k\cdot\a=0\,.
\end{align*}
We note the Galileian invariance upon adding a drift $\c\in\R^3$ and frequency $\omega$,
\begin{align*}
\boldsymbol{v}(t,\boldsymbol{x}) = \psi(t,\boldsymbol{k}\cdot\boldsymbol{x}-\omega t)\boldsymbol{a} + \c\,,
\end{align*}
where the nonlinear term in \eqref{eq 1a} becomes $(\c\cdot\k\frac{\partial\psi}{\partial\xi})\a$, with $\xi=\k\cdot \x-\omega t$ the phase variable of $\psi$, 
which is readily compensated by the time derivative term $-\omega \frac{\partial\psi}{\partial\xi}\a$ for $\omega=\c\cdot\k$; the 
remaining constant term $f\3\times \c$ in \eqref{eq 1a} creates the pressure $p=p_{\c}(x,y):=f(c_2x-c_1y)$, where 
$\c=(c_1,c_2,c_3)^t$.\newline
Such 
vector fields are also always divergence free and therefore solve \eqref{eq 1b}. If $b$ spatially depends on $z$ only, then the 
buoyancy term in \eqref{eq 1a} can be absorbed into the pressure gradient via the primitive $B$ of $b$ with $\frac{d}{dz} B = b$. In 
the barotropic case $v_3\equiv 0$, i.e. $a_3=0$, the nonlinear term in \eqref{eq 1c} also vanishes. What remains are the decoupled 
linear equations with $\xi=\k\cdot\x-\omega t$ and $p=\tilde{p}+B+p_{\c}$,
\begin{align*}
\left(\frac{\partial\psi}{\partial t} - \nu|\k|^2\frac{\partial^2\psi}{\partial \xi^2}\right)\a \ &= \ -f\psi\3\times\boldsymbol{a} - \nabla 
\tilde{p}\\[2mm]
\frac{\partial 
b}{\partial t} - \mu\frac{\partial^2 b}{\partial z^2} \ &= \ 0\,.
\end{align*}
The left hand side of the first equation has the direction $\boldsymbol{a}$, which is orthogonal to $\3\times\a$ on the right hand side, and is 
divergence-free for sufficiently smooth $\psi$, so that both,
left and right hand sides of the first equation, must vanish; this also implies trivial geostrophic balance. The left hand side vanishes for 
$\psi$ solving the heat equation, and the right hand side vanishes for $\boldsymbol{k}=\boldsymbol{e_3}\times\boldsymbol{a}$, 
$\tilde{p}(t,\boldsymbol{x})=-f\Psi(t,\boldsymbol{k}\cdot\boldsymbol{x}-\omega t)$ 
and $\Psi(t,\xi)$ with 
$\frac{\partial\Psi}{\partial\xi}=\psi(t,\xi)$.\newline

In summary, any solutions to the one-dimensional heat equations
\begin{equation}\label{e:heatbasic}
\frac{\partial\psi}{\partial t} = \nu|\boldsymbol{k}|^2\frac{\partial^2 \psi}{\partial \xi^2} \,,\qquad \frac{\partial \tb}{\partial t} = \mu\frac{\partial^2 \tb}{\partial z^2}\,,
\end{equation}
give solutions to the Boussinesq equations \eqref{eq 1}, that we refer to as \emph{horizontal plane flows}, via
\begin{subequations}\label{sol: vanishBoussinesq}
\begin{align}
\boldsymbol{v}(t,\boldsymbol{x}) &= \psi(t,\boldsymbol{k}\cdot\boldsymbol{x}-\omega t)\boldsymbol{a}+\c \qquad \mbox{with}\quad a_3=0 \,,\quad \boldsymbol{k}=\boldsymbol{e_3}\times\boldsymbol{a} \,,\quad \omega=\c\cdot\k\,,\label{sol: vanishBoussinesqa} \\[2mm]
b(t,\boldsymbol{x})&= \tb(t,z) \,,\label{sol: vanishBoussinesqb}\\[2mm]
\begin{split}
p(t,\boldsymbol{x})&=-f\Psi(t,\boldsymbol{k}\cdot\boldsymbol{x}-\omega t)+B(t,z)+p_c(x,y) \\[2mm]
&\mbox{with}\quad \frac{\partial\Psi}{\partial\xi}=\psi(t,\xi)\,,\quad \frac{\partial B}{\partial z}=b(t,z)\,.
\end{split}\label{sol: vanishBoussinesqc}
\end{align}
\end{subequations}
Notably, for each $\a,\k\neq 0$ these plane waves form an injection of solutions to \eqref{e:heatbasic} into the solutions \eqref{sol: 
vanishBoussinesq} of the Boussinesq equations. 
In particular, due to their linear nature, for each nontrivial solution to \eqref{e:heatbasic}, the free choice of prefactors generates a two-dimensional 
set of solutions to \eqref{eq 1}, which is four-dimensional when also counting the free choice of $\a$ with $a_3=0$ (translation and 
Galilean invariance give 4 additional dimensions).

This example illustrates how a plane wave approach reduces the nonlinear equations \eqref{eq 1} effectively to linear equations, and arbitrary superposition of such solutions in the same wave vector direction, but different wave vector lengths $\alpha\k$ with $\alpha\in\R$, are possible. However, clearly general superpositions for different $\k$ create cross terms and thus nonlinear effects.

We emphasize that these horizontal plane flows arise from embedding solutions to the planar Navier Stokes equations given by the first two components of $\v$ as above. Indeed, the velocity field $\boldsymbol{v}$ and buoyancy $b$ do not interact since the velocity field is purely horizontal (no vertical dependence and component),
while the buoyancy is purely vertical (independent of the horizontal directions). While \cite{Majda03} uses the same idea of vanishing nonlinearity and obtains similar solutions, those are without viscosity and Coriolis term. 

\medskip
In the inviscid case $\nu=\mu=0$ the heat equations \eqref{e:heatbasic} just imply time-independence of $\psi$ and $b$, so that these 
shapes of plane waves can be chosen as arbitrary functions of $\xi$.

\medskip
Related, but differently oriented solutions of similar type arise from the parallel flow ansatz $\v(t,\x) = w(t,x,y)\3$, as done in 
\cite{Wang89} for the Navier-Stokes equations. Such a flow is divergence free and the nonlinear term  $(\v\cdot\nabla)\v$ 
vanishes, as does the rotation term $f\3\times\v$. If $b$ is also independent of $z$, then $p=p_0z$ for some $p_0\in\R$ and \eqref{eq 1} 
reduces to the inhomogeneous linear system 
\begin{align*}
\frac{\partial w}{\partial t} =\nu \Delta_h w  + b - p_0\,, \qquad \frac{\partial b}{\partial t} =\mu \Delta_h b +\frac{d\bar{\rho}}{dz} w\,,
\end{align*}
with horizontal Laplacian $\Delta_h$. Solutions can be written as the constant steady state $b=p_0, w=0$ plus superposed Fourier 
modes that decay, and for large scales additionally oscillate in time for stable stratification ($d\bar{\rho}/dz$ constant and negative) 
according the linear dispersion relation; spatially constant modes oscillate in time with frequency $\sqrt{-d\bar{\rho}/dz}$. The 
difference to the horizontal plane flows \eqref{sol: vanishBoussinesq} is not only the vertical velocity direction, but also the coupling of the velocity 
with the buoyancy. However, a superposition with horizontal plane flows \eqref{sol: vanishBoussinesq} is in general not possible, 
since cross terms like
\begin{equation*}
a_1\psi\frac{\partial b}{\partial x} + a_2\psi\frac{\partial b}{\partial y} + w\frac{\partial \tb}{\partial z}
\end{equation*}
remain in \eqref{eq 1c} from $(\v\cdot\nabla)b$. An exception is the single mode parallel flow 
\begin{equation*}
w(t,x,y)=\widehat{w}(t,\k\cdot\x)\,,\quad b(t,x,y)=\hat{b}(t,\k\cdot\x)\,,
\end{equation*}
for $\k=\3\times\a$, whose superposition with horizontal plane flows \eqref{sol: vanishBoussinesq} that have $\tb\equiv0$ yields an explicit solution. 

\medskip
Another related class of solutions are the Kolmogorov flows, as for instance presented in \cite{BalmforthYoung2005} for the non-rotating 
Boussinesq equations. Here a time independent forcing of the single wave mode is implemented in the momentum equation \eqref{eq 1a}, but 
we disregard this for the moment.  Generally, for steady solutions of the form $\v=\cos(kx-mz)\a$, the velocity direction has to be 
$\a=\alpha(m,0,k)^t$ 
and the pressure $p=-\alpha\nu|\a|^2m/k\cdot\sin(kx-mz)$. This yields the buoyancy as 
$b=\alpha\nu|\a|^2(k+m^2/k)\cos(kx-mz)$ 
and the amplitude $\alpha$ is either zero or the stratification is constrained to
\begin{equation*}
\frac{d\bar{\rho}}{dz}=\mu\nu|\a|^4(1+m^2/k^2)\,.
\end{equation*}
Hence, stable stratification (left hand side negative) only allows for the trivial solution $\alpha=0$. Unstable stratification (left hand 
side positive) permits nontrivial solutions with arbitrary amplitudes. In contrast to the parallel flow and the solutions \eqref{sol: 
vanishBoussinesq} these have vertical dependence and vertical velocity component, so that velocity and buoyancy are coupled; 
superposition with the parallel flow or horizontal plane flow \eqref{sol: vanishBoussinesq} do not yield explicit solutions, since their wave vectors and flow directions are not mutually orthogonal, so that cross terms remain from the nonlinear term. 
Note that viscosity and diffusion are required ($\nu\mu\neq 0$) for nontrivial Kolmogorov flow with buoyancy in the presence of stratification. \newline
With forcing of the same form as $\v$ and an amplitude factor $\beta$ in addition to $\alpha$, the solutions are adjusted slightly by the amplitude of the forcing, and nontrivial such solutions must satisfy
\begin{equation*}
\frac{d\bar{\rho}}{dz}=\mu|\a|^2(\nu|\a|^2-\beta)(1+m^2/k^2)\,,
\end{equation*}
while the amplitude of the solution $\alpha$ is still arbitrary. In particular, this also allows for nontrivial solutions in the case of stable 
stratification and shows how forcing can influence the occurrence of steady solutions. In \S\ref{s:forcing} we will study a more general form of forcing and the emerging steady solutions in more detail. 
Note that for $\beta\neq 0$ only $\mu\neq 0$ is required for a nontrivial flow of this type with buoyancy in the presence of stratification.

\medskip
Other related single mode solutions are so-called monochromatic intertia gravity waves (MGW), as described in \cite{Mied76, Drazin77} 
without rotation and in, e.g., \cite{Yau04, Achatz06} for $f\neq 0$. Notably, for $f=0$ these have the same velocity form as the 
Kolmogorov flow, but with phase shifted buoyancy, and these have nonzero temporal frequency. 
MGW exist in the inviscid case $\mu,\nu=0$ and are travelling waves with nonzero 
velocity, depending on a phase variable $\xi = kx-mz-\omega t$ with nonzero squared frequency $\omega^2= (-\frac{d\bar{\rho}}{dz} k^2 + f^2 m^2)/(k^2+m^2)$.
Again, the nonlinear terms vanish, but now the time-derivates compensate the linear terms together with the pressure in 
case $f^2\neq -\frac{d\bar{\rho}}{dz}$. Specifically, the velocity is $\v=\alpha(\cos(\xi)\a + \frac{mf}{\omega}\sin(\xi) 
\boldsymbol{e}_2)$, $\a=(m,0,k)^t$,
$\alpha\in\R$, the pressure $p = \alpha m\frac{\omega^2-f^2}{k\omega}\cos(\xi)$, and the buoyancy $b=-\alpha\frac{k}{\omega}\frac{d\bar{\rho}}{dz}\sin(\xi)$. In the presence of viscosity with $\mu=\nu\neq0$ the MGW turn from stationary into exponentially decaying solutions with the factor for each component given by $\exp(-\nu(k^2+m^2)t)$. 
Exponentially decaying MGWs with $\mu\neq\nu$ are not possible due to the coupling of velocity and buoyancy in this case.\newline 
Superpositions with the parallel flow or horizontal plane flow \eqref{sol: vanishBoussinesq} do not give explicit solutions due to remaining cross terms, exactly as in the case with the Kolmogorov flow. However, for $f=0$ and $\nu=\mu\neq0$ the superposition of MGW and Kolmorogov flow with the same flow direction $\a$ do not give cross terms, since wave vectors and flow directions are orthogonal. Thus, linear terms remain and the superposition of both flows is an explicit solution. Note that in this case the Kolmogorov flow is a steady flow, while the MGW is exponentially decaying as mentioned above.

\subsection{Gradient nonlinearity}\label{s:BoussGrad}
Another special case for the nonlinear term occurs when this is a gradient, thus producing pressure gradient only. 
It seems that \cite{Walsh} was one of the first to notice that any divergence free Laplace eigenfunctions $u_1,u_2$ with the same eigenvalue $\lambda$, i.e., the same wave length, generate a solution $(v_1,v_2)$ to the planar Navier-Stokes equations with viscosity $\nu$. The relation is $v_j=\exp(\lambda \nu t)u_j$, $j=1,2$, and the pressure gradient is given by the negative resulting nonlinear term; a simple case are so-called Taylor flows. 
\cite{Beloshapkin89, Majda02} superpose sinusoidal plane waves with the same wave vector length by summation, while \cite{Hui87} uses integrals, and all of these are causing a gradient nonlinear term.
As in the previous section, these immediately yield solutions to the rotating Boussinesq equations in the form of barotropic flows with decoupled bouyancy. In contrast to the horizontal plane flows \eqref{sol: vanishBoussinesq}, superpositions give explicit solutions when using the same wave length, but arbitrary direction.

We reproduce these two-dimensional solutions from a plane wave viewpoint and additionally provide explicit pressure. For that, analogous to the previous approach, we consider the barotropic case. This means, that the velocity $\v$ exists on the horizontal plane and $b$ spatially depends on $z$ only. For better readability we define here for any $\k\in{\R^2}$ the corresponding three-dimensional vectors $\kb:=(\k,0)^t\in\R^3$ and $\kbp:=(\k^{\perp},0)^t\in\R^3$. 

We start with a simple form of these explicit solutions, which are, like the more general form, superposed generalized Beltrami flows. We use the same approach as for the horizontal plane flow \eqref{sol: vanishBoussinesq}. In order to create a nonzero nonlinear term, we superpose these plane flows, e.g., for two we have
\begin{align*}
\boldsymbol{v}(t,\boldsymbol{x}) = \psi_1(t,\kb_1\cdot\x)\kbp_1+\psi_2(t,\kb_2\cdot\x)\kbp_2\,,
\end{align*}
with, at least for now, arbitrary wave vectors $\k_1, \k_2 \in\R^2$ and wave shapes $\psi_1,\psi_2\in C^{2}(\mathbb{R}_{\geq 0}\times\mathbb{R},\R)$. The nonlinear advective term with this velocity is 
\begin{align*}
(\v\cdot\nabla)\v=\frac{\p\psi_1}{\p\xi}\psi_2\Bigl((\k_1\cdot\k_2)\kb_1-|\k_1|^2\kb_2\Bigr)+\psi_1\frac{\p\psi_2}{\p\xi}\Bigl(-|\k_2|^2\kb_1+(\k_1\cdot\k_2)\kb_2\Bigr)\,.
\end{align*}
If the wave vectors point in the same direction, then the nonlinear term vanishes and we get the horizontal plane flow \eqref{sol: vanishBoussinesq}. If they do not point in the same direction, then, in order to have a gradient nonlinear term, it turns out that the wave vectors and the wave shapes have to satisfy
\begin{align*}
|\k_1|=|\k_2|\,,\quad \psi_i(t,\xi)=\alpha_i(t)\sin(\xi+\delta_i(t))\qquad\mbox{for}\quad i=1,2\,,
\end{align*}
with (for now) arbitrary amplitudes and phase shifts $\alpha_i, \delta_i \in C^1(\mathbb{R}_{\geq 0},\R)$, for $i=1,2$.
For bounded spatial domains with suitable boundary conditions it is also possible to have wave shapes of the form $\alpha(t)e^{\xi}$, but we will not discuss these further since we consider the whole space $\R^3$.
For brevity and readability, in the following $\psi_i$ denotes the sinusoidal wave shape.
With the condition for the wave vectors above the nonlinear term has the gradient structure
\begin{align*}
(\v\cdot\nabla)\v = \nabla h\,,\quad h(t,\x)=(\k_1\cdot\k_2)\psi_1\psi_2+|\k_1|^2\frac{\p\psi_1}{\p\xi}\frac{\p\psi_2}{\p\xi}\,.
\end{align*}
Hence, the pressure can cancel the nonlinear term in \eqref{eq 1a}.\newline
Since the velocity is also divergence free, we can solve the whole three-dimensional system of equations \eqref{eq 1}. The remaining linear system of equations is
\begin{gather*}
\frac{\partial\boldsymbol{v}}{\partial t} - \nu\Delta\boldsymbol{v}= -f\boldsymbol{e_3}\times\boldsymbol{v} - \nabla \tilde{p} \,,\qquad \frac{\partial\tb}{\partial t} - \mu\frac{\partial^2\tb}{\partial z^2} = 0 \,,\\[2mm]
\mbox{where}\quad p(t,\boldsymbol{x})=\tilde{p}(t,\boldsymbol{x})-h(t,\x)+B(t,z) \quad\mbox{with}\quad \frac{\partial B}{\partial z}=\tb(t,z) \,.
\end{gather*}
As before, the left and right hand side of the first equation each have to be zero, but $\nabla h\neq 0$ means there is no geostropic balance on the nonlinear level.
In summary, we find solutions of \eqref{eq 1}, that we refer to as \emph{interacting horizontal plane flows}, given by
\begin{subequations}\label{sol: gradBoussinesq}
\begin{align}
\boldsymbol{v}(t,\boldsymbol{x}) &= \psi_1\left(t,\kb_1\cdot\x\right)\kbp_1+\psi_2\left(t,\kb_2\cdot\x\right)\kbp_2\,,\label{sol: gradBoussinesqa}\\[2mm]
b(t,\boldsymbol{x}) &= \tb(t,z) \,, \label{sol: gradBoussinesqb}\\[2mm]
p(t,\boldsymbol{x}) &= -(\k_1\cdot\k_2)\psi_1\psi_2-|\k_1|^2\frac{\p\psi_1}{\p\xi}\frac{\p\psi_2}{\p\xi}-f\frac{\p(\psi_1+\psi_2)}{\p\xi}+B(t,z)\,,\\[2mm]
\begin{split}
|\k_1|&=|\k_2|\,,\quad \psi_i(t,\xi)=\beta_ie^{-\nu|\k_1|^2t}\sin(\xi+\delta_i)\,,\quad \beta_i,\delta_i\in\R\qquad\mbox{for}\quad i=1,2\,,\\[2mm]
\frac{\partial\tb}{\partial t}&= \mu\frac{\partial^2\tb}{\partial z^2}\,,\quad\frac{\partial B}{\partial z}=\tb(t,z)\,.
\end{split}\label{sol: gradBoussinesqd}
\end{align}
\end{subequations}

\medskip
It is also possible to have superpositions with arbitrary many plane waves, as long as the length of the wave vectors are the same, as done, e.g., for the planar Navier-Stokes equations in \cite{Majda02}. The flow then has the form
\begin{align*}
\v(t,\x)=e^{-\nu a^2t}\sum_{i=1}^{N}\psi_i(\kb_i\cdot\x)\kbp_i \qquad\mbox{with}\quad |\k_i|=a\,,\quad \psi_i(\xi)=\beta_i\sin(\xi+\delta_i)\,,\quad 1\leq i\leq N\,,
\end{align*}
as well as arbitrary constants $\beta_i, \delta_i\in\R$ for $1\leq i\leq N$, and the pressure is
\begin{align*}
p(t,\x)=-e^{-2\nu a^2t}\sum_{i=1}^N\sum_{j=i+1}^N\Bigl((\k_i\cdot\k_j)\psi_i\psi_j+a^2\frac{\p\psi_i}{\p\xi}\frac{\p\psi_j}{\p\xi}\Bigr)\,,
\end{align*}
in order to compensate the resulting nonlinear term.\newline
Analogous to \cite{Hui87}, we may also superpose through an integral over the wave vectors with the same length. For a flow on the horizontal plane
\begin{align*}
\v(t,\x)=e^{-\nu a^2t}\int\displaylimits_{\substack{\ \k\in\R^2 \\ |\k|=a}}\psi_{\k}(\kb\cdot\x)\kbp d\k\qquad\mbox{with}\quad \psi_{\k}=\beta_{\k}\sin(\xi+\delta_{\k})\quad\mbox{for all}\quad |\k|=a\,,
\end{align*}
the pressure, which compensates the nonlinear term, reads
\begin{align*}
p(t,\x)=-a^2e^{-2\nu a^2t}\int_0^{2\pi}\int_{\varphi_1}^{2\pi}\cos(\varphi_1-\varphi_2)\psi_{\varphi_1}\psi_{\varphi_2}+\frac{\p\psi_{\varphi_1}}{\p\xi}\frac{\p\psi_{\varphi_2}}{\p\xi}d\varphi_2d\varphi_1\,,
\end{align*}
with $\k=a(\cos(\varphi),\sin(\varphi))^t$ and $\psi_{\varphi}:=\psi_{\k}$ for all $0\leq\varphi<2\pi$.

\medskip
While for the horizontal plane flows \eqref{sol: vanishBoussinesq} it is possible to superpose arbitrary wave vector lengths of the same wave vector direction,
the interacting horizontal plane flows \eqref{sol: gradBoussinesq} superpose the same wave vector length of arbitrary wave vector directions.
Additionally, the wave modes of the flows \eqref{sol: gradBoussinesq} are interacting with each other, since the nonlinear cross terms do not vanish,
in contrast to the horizontal plane flows \eqref{sol: vanishBoussinesq}, where the cross terms vanish and each wave mode is solving the remaining linear equation on its own,
without influencing the others.
As before, velocity field $\boldsymbol{v}$ and buoyancy $b$ do not influence each other, so that $\v$ as above is also a solution to 
the rotating Navier-Stokes equations and the first two velocity components solve the planar Navier Stokes equations.

\medskip
For better legibility, we do not add a co-moving frame here, even though the Galilean invariance for these solutions can be used in the same way as in \S\ref{s:vanish1};
for arbitrary additional constant velocity term $\c=(c_1,c_2,c_3)^t$, the corresponding non-constant terms attain the wave speeds $\omega_i=\c\cdot\kb_i$ for superposition by summation,
or $\omega_{\k}=\c\cdot\kb$ for superposition by integral, and the linear term $p_c(x,y)=f(c_2x-c_1y)$ is added to the pressure. 
In contrast to the horizontal plane flows \eqref{sol: vanishBoussinesq}, the temporal frequencies here differ from each other, depending on the wave vector direction.

\medskip
As a dimension count for fixed $N$, for each solution $b$ to the heat equation we note that the free amplitudes $\beta_i$ and free scaling of $b$ generate an $N+1$ dimensional set. The wave vectors $\k_i\in\R^2$ of the same length form a one-dimensional set for each $\psi_i$. 
Together with the free shifts $\delta_i$, we have a
$3N+1$ dimensional set of interacting horizontal plane flows \eqref{sol: gradBoussinesq} of the equations \eqref{eq 1} (this includes translation; Galilean invariance gives $3$ more). Clearly, integral superpositions of $\psi_{\k}\kbp$ as above yield an overall infinite dimensional set of solutions.

\medskip
In the inviscid case $\nu=\mu=0$ the heat equation in \eqref{sol: gradBoussinesqd} for $b$ and the form of $\psi_i$ imply time-independence for these wave shapes.
Hence, there is no decay and $b$ is an arbitrary function of $z$.

\section{Explicit solutions in non-rotating fluid models}\label{s:nonrot}
In this section we exploit the ideas of the previous section and illustrate how to obtain explicit superposed plane wave solutions for simpler fluid models. 
The absence of the Coriolis term admits the aforementioned Kolmogorov flow, and in absence of stratification we can set $b\equiv 0$, so 
that the equations become rotationally symmetric. For instance, the parallel flow $\v=w(t,x,y)\3$ can be rotated in any direction and 
yields $\v = w(t,\x)\a$, with arbitrary $\a\in\R^3$ and the invariance $w(t,\x+\alpha \a)=w(t,\x)$ for any $\alpha\in\R$. Likewise, the 
horizontal plane flows \eqref{sol: vanishBoussinesq} and \eqref{sol: gradBoussinesq} with $b\equiv 0$ can be rotated to any direction $\a$.
We focus on the incompressible Navier-Stokes equations on an $n$-dimensional space, $\x\in\R^n$, 
\begin{subequations}\label{eq 2}
\begin{align}
\frac{\partial\boldsymbol{v}}{\partial t}+(\boldsymbol{v}\cdot\nabla)\boldsymbol{v}\ &= \ -\frac{1}{\rho_0}\nabla p +\nu\Delta\boldsymbol{v} \label{eq 2a}\\[2mm]
\nabla\cdot\boldsymbol{v}\ &= \ 0\,,\label{eq 2b}
\end{align}
\end{subequations}
with velocity $\boldsymbol{v}$, pressure $p$, constant density $\rho_0>0$ and kinematic viscosity $\nu>0$;
we also consider the $n$-dimensional Euler equations in \eqref{eq 2} by setting $\nu=0$.\newline

Before presenting the plane wave based solutions, we want to mention the possibility of creating solutions from solutions of lower dimensional problems by separating directions. For illustration of the latter we split the dimension $n=n_1+n_2$ and let $\v^j=(v_1^j,\ldots, v_{n_j}^j)$ with pressure $p_j$ be solutions for dimensions $n_j$, $j=1,2$, respectively. Then, a solution $\v=(v_1,\ldots, v_n)$ in $n$-dimensions with pressure $p$ can be constructed by taking the velocities
$v_k(x_1,\ldots, x_n) := v^1_k(x_1,\ldots, x_{n_1})$ for $k=1,\ldots, n_1$ and 
$v_{n_1+k}(x_1,\ldots, x_n) := v^2_k(x_{n_1+1},\ldots, x_n)$ for $k=1,\ldots, n_2$. The pressure analogously is defined as 
$p(x_1,\ldots,x_n) := p_1(x_1,\ldots, x_{n_1}) + p_2(x_{n_1+1},\ldots, x_n)$. This is a solution, since the transport nonlinearity decouples the $n_1$ and $n_2$ dimensional parts by the separated dependencies of variables. It is also possible to permute the coordinates or rotate the space.\newline
The solutions in \S\ref{s:vanish2} do not belong to this type of constructed solutions, since the wave vectors and flow directions of these solutions are in different subspaces, which are orthogonal to each other. In contrast, the solutions in \S\ref{s:gradnon} belong to this type of constructed solutions, since the flow is decomposed in two-dimensional subspaces, where in each subspace the solutions are of the form as in \cite{Walsh}. However, generalizing them by superpositions with certain parallel flows yields solutions, which do not satisfy this decomposition of variables as well (see Remark \ref{r:supParallelFlow}).

\subsection{Vanishing nonlinearity: `transverse flows'}\label{s:vanish2}
In general, Eulerian fluid models possess the nonlinear terms caused by the material derivative
\begin{align*}
  \frac{{\rm D}\boldsymbol{\zeta}}{{\rm D}t}=\frac{\partial \boldsymbol{\zeta}}{\partial t} + (\boldsymbol{v}\cdot\nabla)\boldsymbol{\zeta}\,,
\end{align*}
where $\boldsymbol{v}=\boldsymbol{v}(t,\boldsymbol{x})\in\mathbb{R}^n$ is the velocity field on the space 
$\boldsymbol{x}\in\mathbb{R}^n$ 
at time $t\geq 0$, with $n\geq 2$, and 
$\boldsymbol{\zeta}=\boldsymbol{\zeta}(t,\boldsymbol{x})\in\mathbb{R}^m$ 
is any quantity transported by the fluid (e.g. density, 
temperature, salinity or the velocity field $\boldsymbol{v}$ itself) with $m\geq 1$. 
To identify flows for which these nonlinear terms vanish, a simple general idea is to make the gradient of each component of $\bzeta$ orthogonal to the velocity $\v$, which means for all $\x\in\R^n$ and $t\geq0$
\begin{equation*}
\v(t,\x)\cdot\nabla\zeta_i(t,\x)=0 \qquad\mbox{for any}\quad 1\leq i \leq m \,.
\end{equation*}
If the velocity always points in one direction, as for the parallel flow $\v(t,\x)=w(t,x,y)\3$ discussed in \S\ref{s:vanish1}, in order for $(\v\cdot\nabla)\bzeta=0$ the quantity $\bzeta$ should spatially depend only on the subspace perpendicular to the velocity direction. The horizontal plane flows \eqref{sol: vanishBoussinesq}, which we will generalize here, depend only on a one dimensional subspace orthogonal to the velocity direction.\newline

Generalising from the approaches of \S\ref{s:bouss} we investigate the superposed plane waves
\begin{align*}
\boldsymbol{v}(t,\boldsymbol{x}) = \boldsymbol{A}\boldsymbol{\psi}(t,\boldsymbol{x}) + \boldsymbol{c}\,,\quad \bzeta(t,\boldsymbol{x}) = \sum_{j=1}^M\boldsymbol{\phi}_{j}(t,\boldsymbol{k}_{j}\cdot\boldsymbol{x}+\omega_j t)\,,
\end{align*}
where $M, N\geq 1$ with $N< n$ and  for $\v$ the flow shapes $\boldsymbol{\psi}\in C^1(\mathbb{R}_{\geq 0}\times\mathbb{R}^n, \R^N)$,
flow directions  $\boldsymbol{A}\in\mathbb{R}^{n\times N}$ and constant drift $\boldsymbol{c}\in\mathbb{R}^n$. The quantity $\bzeta$ contains the wave shapes $\boldsymbol{\phi}_j\in C^1(\mathbb{R}_{\geq 0}\times\mathbb{R},\R^m)$, wave vectors $\boldsymbol{k}_j\in\mathbb{R}^n$, and temporal frequencies $\omega_j\in\mathbb{R}$, where $1\leq j\leq M$. Analogous to before, we presuppose that each wave vector $\boldsymbol{k}_j$ is orthogonal to every column $\boldsymbol{a}_i$ of the matrix $\boldsymbol{A}=[\boldsymbol{a}_1, \ldots, \boldsymbol{a}_N]$, i.e.
\begin{align*}
\boldsymbol{a}_i\cdot\boldsymbol{k}_j=0 \qquad\mbox{for all}\quad 1\leq i\leq N \quad\mbox{and}\quad 1\leq j\leq M\,,
\end{align*}
so that the nonlinear terms caused by the material derivative become
\begin{align*}
(\boldsymbol{v}\cdot\nabla)\boldsymbol{\zeta} = \sum_{j=1}^M(\boldsymbol{c}\cdot\boldsymbol{k}_j)\frac{\partial\boldsymbol{\phi}_j}{\partial\xi}\,,
\end{align*}
where $\xi$ refers to the phase variables of $\boldsymbol{\phi}_j$, respectively.
Notably, the only remnant of $\v$ is the constant $\c$ (the zeroth Fourier mode), so that the nonlinear terms of the material derivative are effectively linear, as a linear combination of the $\xi$-derivatives of the travelling waves of $\boldsymbol{\zeta}$ with the prefactors $\boldsymbol{c}\cdot\boldsymbol{k}_j$. 

When the material derivative acts on the velocity field $\boldsymbol{v}$ itself, we choose all $\bpsi_i$, $i=1,\ldots,N$, in the form of $\bzeta$, so that 
\begin{equation}\label{eq 3}
\boldsymbol{v}(t,\boldsymbol{x}) = \sum_{i=1}^N\boldsymbol{a}_i\sum_{j=1}^{M_i}\psi_{i,j}(t,\boldsymbol{k}_{i,j}\cdot\boldsymbol{x}-\omega_{i,j} t) + \boldsymbol{c}\,,
\end{equation}
where $1\leq N<n$, $M_i\geq1$ and $\boldsymbol{a}_i,\boldsymbol{c}\in\mathbb{R}^n$, $\psi_{i,j}\in C^1(\mathbb{R}_{\geq 0}\times\mathbb{R},\R)$,
$\boldsymbol{k}_{i,j}\in\mathbb{R}^n$, $\omega_{i,j}\in\mathbb{R}$ for all $1\leq j\leq M_i$ and $1\leq i\leq N$. If the wave vectors satisfy the orthogonality conditions
\begin{equation}\label{eq 4}
\boldsymbol{a}_i\cdot\boldsymbol{k}_{j,\ell}=0 \qquad\mbox{for all}\quad 1\leq \ell\leq M_j \quad\mbox{and}\quad 1\leq i,j\leq N\,,
\end{equation}
then the nonlinear term from the material derivative becomes
\begin{equation}\label{eq 5}
(\boldsymbol{v}\cdot\nabla)\boldsymbol{v} = \sum_{i=1}^N\boldsymbol{a}_i\sum_{j=1}^{M_i}(\boldsymbol{c}\cdot\boldsymbol{k}_{i,j})\frac{\partial\psi_{i,j}}{\partial\xi}\,,
\end{equation}
where $\xi$ is the second variable of each $\psi_{i,j}$.
Such a velocity field \eqref{eq 3} is again divergence free due to the orthogonality \eqref{eq 4} and thus suitable for incompressible models. In particular, \eqref{eq 2b} is satisfied.\newline

For simplicity, we presuppose that $p(t,\x)=\tilde{p}(t)$, since the velocity of the form \eqref{eq 3} and \eqref{eq 4} does not produce any gradient terms.
Assuming $\psi_{i,j}\in C^2(\mathbb{R}_{\geq 0}\times\mathbb{R},\R)$ for any wave shape, and substituting the velocity field \eqref{eq 3} with its condition \eqref{eq 4} into the momentum equations of the Navier-Stokes equations \eqref{eq 2a} yields, using \eqref{eq 5}, the linear equations
\begin{align*}
\sum_{i=1}^N\boldsymbol{a}_i\sum_{j=1}^{M_i}\Bigl(\frac{\partial\psi_{i,j}}{\partial t} - \omega_{i,j}\frac{\partial\psi_{i,j}}{\partial\xi} + (\boldsymbol{c}\cdot\boldsymbol{k}_{i,j})\frac{\partial\psi_{i,j}}{\partial\xi}\Bigr) = \sum_{i=1}^N\boldsymbol{a}_i\sum_{j=1}^{M_i}\nu|\boldsymbol{k}_{i,j}|^2\frac{\partial^2\psi_{i,j}}{\partial\xi^2}\,.
\end{align*}
We solve this equation for each index by choosing the frequency and solving the heat equation
\begin{equation}\label{e:heatpsi}
\frac{\partial\psi_{i,j}}{\partial t} = \nu|\boldsymbol{k}_{i,j}|^2\frac{\partial^2\psi_{i,j}}{\partial\xi^2} \,,\quad \omega_{i,j}=\boldsymbol{c}\cdot\boldsymbol{k}_{i,j}\qquad\mbox{for each defined}\quad (i,j)\,.
\end{equation}
By the Galilean invariance (here without effect on the pressure), we may a priori set $\c=0$ when choosing coordinates according to $\omega_{i,j}=\boldsymbol{c}\cdot\boldsymbol{k}_{i,j}$; we note this may be violated in presence of forcing as in \S\ref{s:forcing}.
As in the previous section, each set of solutions to the heat equations \eqref{e:heatpsi} generates an explicit solution to the incompressible Navier-Stokes equations \eqref{eq 2}, which can be verified by a straightforward computation.
We refer to these solutions as \emph{transverse flows}. As already mentioned, these explicit solutions are not constructed by lower-dimensional solutions as shown above, since wave vectors and flow directions belong to different subspaces due to the orthogonality condition \eqref{eq 4}.\newline

For the transverse flow \eqref{eq 3} in the three-dimensional space ($n=3$) there are two cases: We can choose $N=1$, which means that the flow has one direction $\a_1$ and arbitrarily many wave vectors $\k_{1,j}$, lying in the plane orthogonal to the direction $\a_1$. This is related to the parallel flow, cf.\ \cite{Wang89}, but in plane wave form, e.g., $\a_1=\3$ and $\k_{1,j}$ in the horizontal plane. For the case $N=2$ we have two flow directions $\a_1$ and $\a_2$, but only one direction for the wave vectors $\k_{1,j}$ and $\k_{2,j}$, which is orthogonal to the plane spanned by $\a_1$ and $\a_2$. This is similar to a two-dimensional flow with cross flow, cf.\ \cite{Weinbaum67}, but also in a plane wave form. For instance, if the wave vector is in the horizontal plane $\k_{1,j}=\k_{2,j}=(\k,0)^t\in\R^3$, then with $\a_1=(\k^{\perp},0)^t$ we have the velocity of the horizontal plane flow \eqref{sol: vanishBoussinesq} as a purely horizontal flow, superposed with a parallel flow for $\a_2=\3$ as the cross flow component.\newline

In order to make the constraints more clear, we count dimensions of this set of solutions for given $M_i$ (see Remark \ref{r:dense} below for the infinite dimensional case). The dimension from initial data for \eqref{e:heatpsi} that are linearly independent is $m_N:=M_1+\ldots + M_N$ and if these are also linearly independent with respect to scaling for each $1\leq i\leq N$, 
\[
a_1\psi_{i,1}(0,b_1 \xi)+ \ldots + a_{M_i}\psi_{i,M_i}(0,b_{M_i} \xi) \neq 0\,,
\]
not all $a_j,b_j=0$, this gives linearly independent summands in \eqref{eq 3}.
Additionally, we get $N$ dimensions from linearly independent $\a_i$, and different wave vectors also generate linear independence. Hence, admissible wave vectors $\k_{i,j}$ satisfying \eqref{eq 4} can be independently selected from the orthogonal complement of $\mathrm{span}(\a_1,\ldots,\a_N)$ in $\R^n$, which has dimension $(n-N)$ and thus contributing $(n-N)m_N$ dimensions to a total of $N+(n-N+1)m_N$
(Galileian invariance and translation in space give another $2n-N$ dimensions).\newline
As mentioned for the three dimensional problem ($n=3$), we have the two cases $N=1,2$. For $N=1$ the dimension is $1+(3-1+1)m_1=1+3m_1$, while for $N=2$ it is $2+(3-2+1)m_2=2+2m_2$ (and additionally $6-N$ by Galileian invariance and translation in space). For the two dimensional problem we can only choose one flow direction ($N=1$), which yields the dimension $1+(2-1+1)m_1=1+2m_1$ (and additional 3 for Galileian invariance and translation in space).\newline

\begin{remark}\label{r:combine}
The representation in the form of \eqref{eq 3} can be reduced by combining summands of fixed $i$, for which the wave vectors lie on the same ray, e.g. $\k_{i,j} = r\k_{i,\ell}$, $r\in\R$.
However, the combined function does not necessarily solve a heat equation.
\end{remark}

\begin{remark}\label{r:dense}
The space of solutions is infinite dimensional, since the numbers of summands $M_i$ are arbitrary. As before, the summation over $j$ can in fact be replaced by an integral over all admissible wave vectors orthogonal to all flow directions $\a_i$. 
\end{remark}

\medskip
The eigenmode solutions of heat equations yield the simplest solutions from the above class with arbitrary $\alpha_{i,j},\delta_{i,j}\in\mathbb{R}$ for any $1\leq j\leq M_i$ and $1\leq i\leq N$ as
\begin{align*}
\boldsymbol{v}(t,\boldsymbol{x}) = \sum_{i=1}^N\boldsymbol{a}_i\sum_{j=1}^{M_i}\alpha_{i,j}e^{-\nu|\boldsymbol{k}_{i,j}|^2t}\sin(\boldsymbol{k}_{i,j}\cdot\boldsymbol{x}-\omega_{i,j} t + \delta_{i,j}) + \boldsymbol{c}\,.
\end{align*}

\bigskip
For the Euler equations we may proceed in the same way, setting $\nu=0$, so that \eqref{e:heatpsi} becomes $\frac{\partial\psi_{i,j}}{\partial t} = 0$. Hence, the shapes of the travelling wave components of the solutions can be arbitrary functions of $\xi$ and do not explicitly depend on $t$.\newline

As indicated in \S\ref{s:vanish1}, for this class of solutions, adding a Coriolis term poses additional constraints on admissible direction and wave vectors. These are satisfied, e.g., for zero vertical component and $\k=\3\times\a$, so the pressure gradient balances the Coriolis term. The parallel flow, a horizontally dependent vertical flow $\v(t,\x)=w(t,x,y)\3$ \citep{Wang89}, has even a vanishing Coriolis term.

\subsection{Gradient nonlinearity: `interacting transverse flows'}\label{s:gradnon}
As for the Boussinesq equations in \S\ref{s:BoussGrad}, another approach to deal with the transport nonlinearity is to involve the pressure gradient. For illustration, we first present an approach which, in general, does not provide divergence free solutions, but is nevertheless instructive. We start with a velocity field that is a superposition of travelling waves, where the direction of each is  its own wave vector, i.e.
\begin{align*}
\boldsymbol{v}(t,\boldsymbol{x}) =  \sum_{i=1}^N \psi_i(t,\boldsymbol{k}_i\cdot\boldsymbol{x}-\omega_it)\boldsymbol{k}_i+ \boldsymbol{c}\,,
\end{align*}
for any $N\geq 1$ and arbitrary $\psi_i\in C^1(\mathbb{R}_{\geq0}\times\mathbb{R},\R)$ as well as $\boldsymbol{k}_i$, $\boldsymbol{c}\in\mathbb{R}^n$ and $\omega_i\in\mathbb{R}$ for any $1\leq i\leq N$. Then the nonlinear term from the material derivative becomes
\begin{align*}
(\boldsymbol{v}\cdot\nabla)\boldsymbol{v} = \sum_{i=1}^N (\boldsymbol{c}\cdot\boldsymbol{k}_i)\frac{\partial\psi_i}{\partial\xi}\boldsymbol{k}_i + \nabla h\,,
\end{align*}
where $\xi$ is the second variable of each $\psi_i$ and the scalar valued function $h$ is defined as
\begin{align*}
h(t,\boldsymbol{x}) = \frac{1}{2}\sum_{i=1}^N\sum_{j=1}^N (\boldsymbol{k}_i\cdot \boldsymbol{k}_j)\psi_i(t,\boldsymbol{k}_i\cdot\boldsymbol{x}-\omega_it)\psi_j(t,\boldsymbol{k}_j\cdot\boldsymbol{x}-\omega_jt)\,.
\end{align*}
With this choice of $\boldsymbol{v}$ we just have to define the pressure, e.g., in \eqref{eq 2a}, as $p=-\rho_0h$, and the nonlinear term is removed by the pressure gradient $\nabla p$.
By definition of $\v$, the nonlinear terms generated by travelling waves with linearly dependent wave vectors and directions can all be compensated by the pressure gradient.
However, these functions are in general not divergence free and thus do not solve incompressible fluid equations.

\medskip
In order for the velocity field to be divergence-free, we combine the last idea with that of \S\ref{s:vanish2} and, roughly speaking, choose a composition of travelling waves,
where the wave vector of each travelling wave is orthogonal to its flow direction.
As for the interacting horizontal plane flows \eqref{sol: gradBoussinesq}, the interacting travelling waves will depend on the same two-dimensional linear subspace of $\R^n$,
i.e., the wave vectors and flow directions of the interacting waves belong to the same subspace.
This ensures that the nonlinear term resulting from the interacting travelling waves form a gradient, and additionally we are able to describe the nonlinear term explicitly. 
As already mentioned, the resulting solutions are built from those of the purely planar case, cf.\ \cite{Walsh}, using the decomposition of coordinates and rotation mentioned in the beginning of \S\ref{s:nonrot}. Hence, the planar flows are independently placed in two-dimensional orthogonal subspaces of $\R^n$ -- for $n=3$ this is a single rotated plane. 

Specifically, we define for $1\leq N\leq n/2$ arbitrary two-dimensional linear subspaces $S_i\subset\R^n$, for $i=1,\dots,N$, where each subspace $S_i$ is orthogonal to any other subspace $S_j$, with $i\neq j$. Then we consider a velocity field of the form
\begin{equation}\label{eq 6}
\boldsymbol{v}(t,\boldsymbol{x}) =  \sum_{i=1}^N\sum_{j=1}^{M_i} \psi_{i,j}(t,\k_{i,j}\cdot\x-\omega_{i,j}t)\a_{i,j} +\c\,,
\end{equation}
where $\boldsymbol{c}\in\R^n$, $M_i\geq 1$ for any $1\leq i\leq N$, and $\omega_{i,j}\in\R$ for any $1\leq j\leq M_i\,,\,1\leq i\leq N$, with the conditions 
\begin{align}\label{cond: gradSol}
\k_{i,j}, \a_{i,j}\in S_i\,,\quad \k_{i,j}\cdot\a_{i,j}=0\,,\quad |\k_{i,j}|=\lambda_i\,,\quad \psi_{i,j}(t,\xi)=\alpha_{i,j}(t)\sin(\xi+\delta_{i,j}(t))\,,
\end{align}
where $\lambda_i\in\R$ are arbitrary as well as (at least for now) $\alpha_{i,j}, \delta_{i,j}\in C^1(\R_{\geq0},\R)$. As in \S\ref{s:BoussGrad} it is possible to have wave shapes $\psi_{i,j}$ of the form $\alpha(t)e^{\xi}$ in bounded spatial domains with suitable boundary conditions, which we will not discuss further, since we consider the whole space $\R^n$. Again, for brevity and readability, $\psi_{i,j}$ denotes the sinusoidal wave shape in the following.
With \eqref{eq 6} and \eqref{cond: gradSol} the nonlinear term caused by the material derivative of $\boldsymbol{v}$ becomes
\begin{equation}\label{eq 8}
(\boldsymbol{v}\cdot\nabla)\boldsymbol{v} = (\boldsymbol{c}\cdot\nabla)\boldsymbol{v} +\nabla h\,,
\end{equation}
where the scalar valued function $h$ is defined as
\begin{equation}\label{eq 9}
h(t,\boldsymbol{x}) = \sum_{i=1}^N\sum_{j=1}^{M_i}\sum_{\ell=j+1}^{M_i}\left((\k_{i,j}\cdot\k_{i,\ell})\psi_{i,j}\psi_{i,\ell}+\lambda_i^2\frac{\p\psi_{i,j}}{\p\xi}\frac{\p\psi_{i,\ell}}{\p\xi}\right)\,.
\end{equation}
\newline
Recall that in the approach of \S\ref{s:vanish1} and \S\ref{s:vanish2} the nonlinear term vanishes, because the travelling waves do not interact due to the orthogonality condition \eqref{eq 4}. In the approach here and in \S\ref{s:BoussGrad} the travelling waves do interact, but in such a way, that the arising nonlinear term is a gradient and therefore can be compensated by the pressure. The form of the nonlinear interaction $h$ in \eqref{eq 9} shows, that with the approach \eqref{eq 6} and \eqref{cond: gradSol} the travelling waves within the same two-dimensional linear subspace $S_i$ interact as pairs, but due to the orthogonality there is no interaction with those in the other subspaces $S_j$, $i\neq j$.
These interactions are constrained to the same spatial scale, since the wave vectors of the interacting travelling waves have the same length $\lambda_i$. However, the wave vector lengths $\lambda_i$ can be different for each subspace $S_i$.

Upon substituting the velocity field $\boldsymbol{v}$ in the form \eqref{eq 6} with its conditions \eqref{cond: gradSol} and the pressure $p=\tilde{p}-\rho_0h$, $h$ as in \eqref{eq 9}, into the momentum equations of the Navier-Stokes equations \eqref{eq 2a}, we obtain, using \eqref{eq 8}, the linear equations
\begin{align*}
\frac{\partial\boldsymbol{v}}{\partial t}+(\boldsymbol{c}\cdot\nabla)\boldsymbol{v}&=-\frac{1}{\rho_0}\nabla\tilde{p} +\nu\Delta\boldsymbol{v}\,.
\end{align*}
These can be readily solved for the above superposed travelling wave form of $\v$, which yields \eqref{eq 6} with \eqref{cond: gradSol}, \eqref{eq 9} and for any $1\leq j\leq M_i, \ 1\leq i\leq N$
\begin{equation}\label{cond: gradSol2}
p(t,\x)=\tilde{p}(t)-\rho_0h(t,\x)\,,\quad \alpha_{i,j}(t)=\beta_{i,j}e^{-\nu\lambda_i^2t}\,,\quad \beta_{i,j}, \delta_{i,j} \in\R\,,\quad \omega_{i,j}=c\cdot\k_{i,j}\,,
\end{equation}
with arbitrary $\tilde{p}\in C(\R_{\geq0},\R)$. We refer to these solutions as \emph{interacting transverse flows}.
Again, these $\omega_{i,j}$ also follow from the Galilean invariance.  Remark that $\v$ satisfies \eqref{eq 2b} due to the form \eqref{eq 6} and the orthogonality relations described in \eqref{cond: gradSol}.\newline

We thus obtain a set of solutions of the incompressible Navier-Stokes equations, where the wave shapes are sinusoidal and exponentially decaying in time. Notably, the pressure is decaying as the product of the interacting travelling waves, faster than each of these. With suitable forcing the decay of solutions can be compensated, cf.\ \S\ref{s:forcing}.\newline

For this set of solutions a dimension count is as follows: With the maximum value $N= \lfloor n/2 \rfloor$ we have the dimension $m_N:=M_1+\dots+M_N$ for the amplitudes $\beta_{i,j}$, as well as for the phase shifts $\delta_{i,j}$ (including translation in space). The wave vectors $\k_{i,j}$ are in a plane $S_i$ and of the same length $\lambda_i$, which gives one dimension for each wave shape $\psi_{i,j}$. 
Together, we thus count the dimension $3m_N$ (Galileian invariance with constant vector $\c$ adds $n$), which makes $3m_1$ for $n=3$.\newline

\begin{remark}\label{r:denseGrad}
As before, we actually have an infinitely dimensional solution space. 
The inner sum of \eqref{eq 6} over $j$ can be replaced by a general integral over all wave vectors of the same two-dimensional subspace $S_i$ with the same length $\lambda_i$
\begin{align*}
\v(t,\x)=\sum_{i=1}^Ne^{-\nu\lambda_i^2t}\int\displaylimits_{\substack{\k\in S_i \\ |\k|=\lambda_i}}\psi_{i,\k}(\k\cdot\x-\omega_{i,\k}t)\a_{i,\k} d\k +\c \,,
\end{align*}
with similar conditions as before in order to get a gradient nonlinear term,
\begin{align*}
\psi_{i,\k}(\xi)=\beta_{i,\k}\sin(\xi+\delta_{i,\k})\,,\quad \k,\a_{i,\k}\in S_i\qquad \mbox{with}\quad |\k|=\lambda_i\quad \mbox{and}\quad \k\cdot\a_{i,k}=0\,,\quad 1\leq i\leq N\,.
\end{align*}
In order to explicitly determine the corresponding pressure, we first choose for any subspace $S_i$ two orthogonal unit vectors $\boldsymbol{e_{i,1}}, \boldsymbol{e_{i,2}}\in S_i$,
so that any vector $\k_i\in S_i$ is clearly defined by an angle $0\leq\varphi<2\pi$ and factor $\lambda_i>0$ with $\k_i=\lambda_i(\cos(\varphi)\boldsymbol{e_{i,1}}+\sin(\varphi)\boldsymbol{e_{i,2}})$.
Without loss of generality we assume for the corresponding $\k_i$, that $\a_{i,\k_i}=\lambda_i(-\sin(\varphi)\boldsymbol{e_{i,1}}+\cos(\varphi)\boldsymbol{e_{i,2}})$,
since wave vector and flow direction are in the same two-dimensional subspace and the length of $\a_{i,\k_i}$ can be absorbed into the wave shape $\psi_{i,\k_i}$.
Then the pressure, which is compensating the resulting nonlinear term, reads
\begin{align*}
p(t,\x)=\sum_{i=1}^N-\lambda_i^2e^{-2\nu\lambda_i^2t}\int_0^{2\pi}\int_{\varphi_1}^{2\pi}\cos(\varphi_1-\varphi_2)\psi_{i,\varphi_1}\psi_{i,\varphi_2}+\frac{\p\psi_{i,\varphi_1}}{\p\xi}\frac{\p\psi_{i,\varphi_2}}{\p\xi}d\varphi_2d\varphi_1\,,
\end{align*}
with $\k_i=\lambda_i(\cos(\varphi)\boldsymbol{e_{i,1}}+\sin(\varphi)\boldsymbol{e_{i,2}})$ and $\psi_{i,\varphi}:=\psi_{i,\k_i}$ for all $0\leq\varphi<2\pi$ and $1\leq i\leq N$.
\end{remark}

\begin{remark}\label{r:supParallelFlow}
It is possible to superpose the interacting transverse flows \eqref{eq 6} with certain parallel flows. Considering the case $n=3$ and $S_1=\R^2\times\{0\}$ the superposition is
\begin{align*}
\v(t,\x)=e^{-\nu\lambda^2t}\sum_{j=1}^M\alpha_j\left(\sin(\kb_j\cdot\x-\omega_jt+\delta_j)\kbp_j+\gamma\cos(\kb_j\cdot\x-\omega_jt+\delta_j)\3\right)+\c\,,
\end{align*}
with $M\geq1$, $\alpha_j, \delta_j, \gamma\in\R$, $\c\in\R^3$, $\kb_j\in S_1$ satisfying $|\kb_j|=\lambda$ and $\kbp_j$ as defined in \S\ref{s:BoussGrad}, $\omega_j=\c\cdot\kb_j$ for all $j=1, \dots, M$. Since the resulting nonlinear terms from the additional parallel flow components cancel each other, the pressure remains the same as for the interacting transverse flows \eqref{eq 6}, so $p(t,\x)=\tilde{p}(t)-\rho_0h(t,\x)$ with $h(t,\x)$ as in \eqref{eq 9}. These solutions are also presented in \cite{Chai20} using helical decomposition rather than the plane wave ansatz.\newline
The superposition can also be of integral form analogous to Remark \ref{r:denseGrad} for $n=3$, i.e.,
\begin{align*}
\v(t,\x)=e^{-\nu\lambda^2t}\int\displaylimits_{\substack{\k\in\R^2 \\ |\k|=\lambda}}\alpha_{\k}\left(\sin(\kb\cdot\x-\omega_{\k}t+\delta_{\k})\kbp + \gamma\cos(\kb\cdot\x-\omega_{\k}t+\delta_{\k})\3\right) d\k +\c \,,
\end{align*}
with $\gamma, \alpha_{\k}, \delta_{\k}\in\R$, $\kb=(\k,0)^t\in S_1$ and $\omega_{\k}=\c\cdot\kb$ for all $\k\in\R^2$, $|\k|=\lambda$. The pressure remains the same as in Remark \ref{r:denseGrad} for $N=1$, since the additional nonlinear terms cancel each other. These superposed solutions can be rotated, so that $S_1$ need not be horizontal.\newline
This construction generalises to $n>3$ as above, using interacting transverse flows of the form \eqref{eq 6} or Remark \ref{r:denseGrad}. For this one adds to each travelling wave component on $S_i$ a parallel flow component of the same form and on the same plane, but phase shifted to a cosine and a flow direction in $\R^n\backslash S$, with $S$ the span of all $S_i$. The flow directions of the parallel flows can differ for different subspaces $S_i$, but are the same for any component that lies on the same $S_i$, for each $i=1,\dots,N$. These superposed flows are not of the decomposed variable form from the beginning of \S\ref{s:nonrot}.
\end{remark}

\medskip
For the Euler equations we simply set $\nu=0$ in \eqref{eq 2}, which yields the same form of solutions, but with constant amplitudes $\alpha_{i,j}(t)\equiv\beta_{i,j}$ in the conditions \eqref{cond: gradSol2}, i.e., time independent solutions -- as expected.

\medskip
As already mentioned in \S\ref{s:BoussGrad}, the additional Coriolis term in the rotating Boussinesq equations poses additional constraints on this set of solutions.
These are satisfied, e.g., for purely horizontal wave vectors and flow directions -- more precisely, $\k_{1,j},\a_{1,j}\in\R^3$ with $(\k_{1,j})_3=0$ and $\a_{1,j}=\3\times\k_{1,j}$. In that way, the pressure gradient balances the Coriolis term.

\section{Explicit solutions with adapted forcing}\label{s:forcing}
In this section we consider incompressible fluid equations with certain forcing, which can be treated with the same methods as presented before, yielding explicit solutions of forced nonlinear fluid equations. 
Here we consider an $n$-dimensional incompressible fluid equation
\begin{subequations}\label{eq 10}
\begin{align}
\frac{\partial\boldsymbol{v}}{\partial t}+(\boldsymbol{v}\cdot\nabla)\boldsymbol{v} \ &= \ -\nabla p +L\boldsymbol{v} + \boldsymbol{F}\label{eq 10a}\\[2mm]
\nabla\cdot\boldsymbol{v} \ &= \ 0\,,\label{eq 10b}
\end{align}
\end{subequations}
with linear differential operator in space $L$ and forcing $\boldsymbol{F}=\boldsymbol{F}(t,\boldsymbol{x})$.
Based on vanishing or gradient nonlinear terms, we will generate explicit solutions to \eqref{eq 10} in the presence of suitable forcing as solutions to the linear inhomogeneous equation
\begin{equation}\label{e:inhom}
\frac{\partial\boldsymbol{v}}{\partial t}=L\boldsymbol{v} + \boldsymbol{F}\,.
\end{equation}

In detail, we consider functions $\F, \v$ in the form \eqref{eq 3} or \eqref{eq 6} respectively,
with the corresponding orthogonality or length conditions for the wave vectors and flow directions,
and choose $\boldsymbol{c}=0$.
We note that simpler solutions and forcing in this context has been considered, e.g., in \cite{Meshalkin61, BalmforthYoung2005, Beloshapkin89}. Substitution of velocities in the form \eqref{eq 3} or \eqref{eq 6} with its conditions into \eqref{eq 10a} yields the inhomogeneous linear equation \eqref{e:inhom},
since the pressure cancels the nonlinear term, which may be of gradient form. With this approach \eqref{eq 10b} is also satisfied.\newline
We assume that $L$ generates a semigroup $e^{Lt}$ and that this respects the form \eqref{eq 3} or \eqref{eq 6} respectively, with its corresponding conditions for wave vector directions or lengths, e.g., $L$ may consist of constant coefficient differential operators such as the viscosity term $L=\nu\Delta$, which yields the heat-semigroup. Then the homogeneous solution reads
\begin{align*}
\boldsymbol{v}_{\rm hom}(t,\boldsymbol{x})=\left(e^{Lt}\boldsymbol{v}_0\right)(\boldsymbol{x})\,, \quad t\geq 0\,,
\end{align*}
and has the form \eqref{eq 3} or \eqref{eq 6} respectively, with its corresponding orthogonality or wave vector length conditions, if the initial condition $\boldsymbol{v}_0$ does. We then solve the inhomogeneous equation with the particular solution
\begin{align*}
\boldsymbol{v}_{\rm p}(t,\boldsymbol{x})=\int_{0}^te^{L(t-s)}\boldsymbol{F}(s,\boldsymbol{x})ds\,,
\end{align*}
assuming that the forcing $\boldsymbol{F}$ is such that this integral exists. Then \eqref{e:inhom} with initial condition $\v_0$ possesses the solution
\begin{equation}\label{e:inh_sol}
\boldsymbol{v}(t,\boldsymbol{x}) = \boldsymbol{v}_{\rm hom}(t,\boldsymbol{x}) + \boldsymbol{v}_{\rm p}(t,\boldsymbol{x})= \left(e^{Lt}\boldsymbol{v}_0\right)(\boldsymbol{x}) + \int_{0}^te^{L(t-s)}\boldsymbol{F}(s,\boldsymbol{x})ds\,.
\end{equation}
The key observation is that this is indeed a solution to  \eqref{eq 10a} whenever $\boldsymbol{F}$ has the form \eqref{eq 3} or \eqref{eq 6} respectively, and satisfies the corresponding orthogonality and wave vector length conditions, together with $\v_0$. In particular, this ensures that $\v$ is divergence free, i.e., \eqref{eq 10b} is also satisfied. In other words, the dynamics for initial data and forcing under these constraints is linear.

As usual for inhomogeneous linear equations, if $\F$ is in addition time independent and lies in the range of $L$, then any preimage $\v_s = L^{-1} \F$ is a steady state and $\w = \v-\v_s$ solves the homogeneous equation $\w_t = L\w$, for any solution $\v$ of \eqref{e:inhom}. Hence, $\w = e^{Lt}\w_0$ with $\w_0 = \v_0 - \v_s$, so that \eqref{e:inh_sol} can be cast as $\v = e^{Lt}(\v_0-\v_s) + \v_s$. The stability of any such steady state $\v_s$ \emph{within} the linear subspace of \eqref{eq 3} or \eqref{eq 6}, respectively, is thus directly determined by spectral properties of $L$. 
The simplest case of forcing is an eigenfunction $\F$ of $L$ with nonzero real eigenvalue, $L\F=\lambda^{-1} \F$ with $\lambda\in\R\backslash\{0\}$, for which $\v_s=-\lambda \F$ gives a stationary solution to \eqref{e:inhom} and to \eqref{eq 10}.\newline

We explain the last results with the simple case $L=\nu\Delta$ and a time-independent forcing $\F(\x)$. A well known example for explicit solutions in this case is the aforementioned Kolmogorov flow, as for instance studied regarding stability properties in a non-rotating Boussinesq setting in \cite{BalmforthYoung2002,BalmforthYoung2005}. More generally, let $\a_F,\k_F\in\R^n$ and consider
\begin{align*}
L=\nu\Delta\,, \quad\F(\x)=\int_0^{\infty}\alpha_F(\xi)\sin(\xi\k_F\cdot\x)d\xi\;\a_F \qquad\mbox{with}\quad|\k_F|=1\,, \quad\a_F\cdot\k_F=0\,,
\end{align*}
with integrable density function $\alpha_F:\R_{\geq0}\rightarrow\R$, so that
\begin{align*}
\int_0^1|\alpha_F(\xi)|d\xi+\int_1^{\infty}\xi^2|\alpha_F(\xi)|d\xi<\infty\,.
\end{align*}
Let $\v_0$ be an initial condition of the analogous form with $\a_0,\k_0\in\R^n$ and 
\begin{align*}
\v_0(\x)=\int_0^{\infty}\alpha_0(\xi)\sin(\xi\k_0\cdot\x)d\xi\;\a_0\qquad\mbox{with}\quad
|\k_0|=1 \,,\quad \a_0\cdot\k_0=\a_0\cdot\k_F=\a_F\cdot\k_0=0 \,,
\end{align*}
and integrable density function $\alpha_0:\R_{\geq0}\rightarrow\R$ satisfying the same bound as $\alpha_F$.
Then the corresponding solution of \eqref{e:inhom}, given by \eqref{e:inh_sol}, can be written as
\begin{align*}
\v(t,\x)=&\int_0^{\infty} e^{-\nu\xi^2t}\alpha_0(\xi)\sin(\xi\k_0\cdot\x)d\xi\;\a_0
+\int_0^t\int_0^{\infty}e^{-\nu\xi^2(t-s)}\alpha_F(\xi)\sin(\xi\k_F\cdot\x)d\xi ds\; \a_F\,.
\end{align*}
Alternatively, using the steady particular solution of \eqref{e:inhom} given by 
\begin{align*}
\v_s(\x)= \int_0^{\infty}\frac{\alpha_F(\xi)}{\nu\xi^2}\sin(\xi\k_F\cdot\x)d\xi\; \a_F\,,
\end{align*}
additionally assuming sufficiently quick decay of the density $\alpha_F$ near $\xi=0$, the solution above can be reformulated as 
\begin{align*}
\v(t,\x)= \v_s(\x) + \int_0^{\infty} e^{-\nu\xi^2t}\left(\alpha_0(\xi)\sin(\xi\k_0\cdot\x)\a_0- \frac{\alpha_F(\xi)}{\nu\xi^2}\sin(\xi\k_F\cdot\x)\a_F\right) d\xi\,.
\end{align*}
It is an explicit solution of \eqref{e:inhom} and \eqref{eq 10}, since the nonlinear term vanishes with this choice of wave vectors, as well as flow and forcing directions. The last equation in particular shows, that for this $\F$ the unique (up to a constant) bounded equilibrium state $\v_s(\x)$ of \eqref{e:inhom} is asymptotically stable within the subspace of linear dynamics for \eqref{eq 10}.\newline

The forcing, as well as the velocity, are constructed in line with the transverse flow solutions in \S\ref{s:vanish2}, therefore they are superposed by modes, which have the same direction, but different eigenvalues. According to the properties of the solutions in \S\ref{s:vanish2}, it is also possible to use forcing with different directions, as long as the orthogonality conditions are satisfied. As mentioned in Remark \ref{r:dense}, integration over the different directions and eigenvalues instead of summation is possible.

In the similar way as shown above, one can also construct forcing and velocity in the form of interacting transverse flows, like in \S\ref{s:gradnon}.
Instead of integrating over different eigenvalues and the same direction, one integrates over the same eigenvalue and different directions on a plane (or several planes for $n>3$).
The nonlinear term in \eqref{eq 10} is a gradient and compensated by the pressure in that way, so that the linear inhomogeneous equation \eqref{e:inhom} remains again.\newline

\begin{remark}\label{r:forcedGrad}
In \eqref{e:inhom} we assumed the pressure gradient precisely cancels the nonlinear term. However, we may include a remaining pressure gradient $\nabla\tilde{p}$, if it has the plane wave form. In this case the forcing term takes the form $\F = \tilde{\F} -\nabla\tilde{p}$, where in the simplest case we require $\tilde{p}:=\delta\Psi(\k_F\cdot\x)$ with $\frac{\p\Psi}{\p\xi}=\psi$ and $\tilde{\F}:=\psi(\k_F\cdot\x)\a_\delta$ with $\a_{\delta}:=\a_F+\delta\k_F$. 
Then $\F=\psi(\k_F\cdot\x)\a_F$, so that we have the same situation as above and can proceed in the same way.
\end{remark}

\medskip
We remark that \cite{Chae96} find travelling wave-like solutions in an $n$-dimensional space with linear forcing term,
that -- in contrast -- is a time dependent factor of the velocity.

We can also use the results of this section in the Boussinesq equations with forcing in the momentum equations \eqref{eq 1a}.
As in the horizontal plane flows \eqref{sol: vanishBoussinesq} and \eqref{sol: gradBoussinesq} the velocity and the buoyancy are decoupled by a suitable choice of wave vector and flow direction, so that a two dimensional version of \eqref{eq 10} with Coriolis term remains. The Coriolis term can also be compensated by the pressure gradient as for \eqref{sol: vanishBoussinesq} and \eqref{sol: gradBoussinesq}. A case of solutions of forced Boussinesq equations with coupled buoyancy and velocity is the aforementioned Kolmogorov flow.

\section{Discussion}\label{s:discuss}
In this paper we have presented various kinds of explicit solutions in nonlinear fluid models, which fall into the class of generalized Beltrami flows. We have also investigated possible superpositions, forcings and their linear dynamics within the nonlinear equations.

In \S\ref{s:bouss} we have started with a simple example and an overview of related explicit horizontal plane flows in the rotating Boussinesq equations, distinguished by vanishing or gradient form of the transport nonlinearity. These are based on the orthogonality of wave vector and plane wave directions, and the velocity field also solves the planar Navier Stokes or Euler equations, since purely horizontal velocity and purely vertical buoyancy are decoupled.
In all cases a system of linear equations remains, after compensating the possibly resulting gradient nonlinear term by the pressure gradient.
The same approach for vanishing nonlinearity has been used for Kolmogorov flow \citep{BalmforthYoung2005} and monochromatic inertia gravity waves (MGW) \citep{Mied76,Yau04},
but with coupled velocity and buoyancy, and similarly in parallel flow \citep{Wang89}; in the discussion of these types of flows we have identified an explicit superposition of MGW and Kolmogorow flow, that we have not found elsewhere, cf.\ \S\ref{s:bouss}.

In \S\ref{s:nonrot} we have then generalized these explicit (interacting) horizontal plane flows to the simpler Navier-Stokes and Euler equations as (interacting) transverse flows.
There we have shown, that for explicit solutions with vanishing nonlinear terms it is possible to make superpositions with arbitrary wave lengths, as long as the orthogonality conditions of the wave vectors and plane wave directions are satisfied.
In three dimensions they fall into the class of parallel flow \citep{Wang89} or, more generally, the two-dimensional flow with cross flow \citep{Weinbaum67}.
For explicit solutions with gradient nonlinear terms it is the other way round. The plane waves are superposed for arbitrary wave vector directions on a plane, but with the same wave length.
In the three-dimensional case these are two-dimensional explicit solutions as in \cite{Walsh}, which can be rotated in space.
In higher dimensions one can separate the space in two-dimensional orthogonal subspaces, in which these two-dimensional flows exist.
The explicit solutions in different subspaces do not interact with each other and can have different wave lengths.
These solutions can be generalized by superpositions with certain parallel flows, so that the spatial dependencies are still restricted to the two-dimensional subspaces, but not the flow directions (see Remark \ref{r:supParallelFlow} or for the 3D case \cite{Chai20}).
In \S\ref{s:forcing} we have then added suitable forcing into the momentum equations, which admits explicit steady state solutions of the type as shown in \S\ref{s:nonrot}.
Herein we have allowed the viscous term to be replaced by a more general linear term and we have solved the remaining inhomogeneous linear problem assuming a semi-group for the homogeneous part.
This approach can be readily extended to additional wave vectors or flow directions, and also to the Boussinesq equations.

\medskip
We next discuss the relation to previous results in the literature in some more detail, organised by equation and dimension. For the 2D Euler equations, \cite{Majda02} characterise in Prop.\ 2.2 solutions via a nonlinear Poisson equation for the streamfunction. \cite{Weinbaum67, Wang89, Wang90} analogously consider solutions to 2D Navier-Stokes equations. \cite{Hui87} even adds an $y$-dependent linear term, causing the solutions to exist on a special co-moving frame. Our solutions intersect these sets, but also arise in higher dimensions, and generally do not satisfy the conditions of stream function and vorticity, in particular $\Delta\psi(t,x)=F(\psi(t,x))$ from \cite{Majda02}. These conditions do not allow the superposition of explicit plane wave solutions with different wave length, which is indeed possible in our solutions; notably the `only if' in Prop.\ 2.2 in \cite{Majda02} is too strong.\newline
Nevertheless, as mentioned, \cite{Hui87, Beloshapkin89} and \cite{Majda02} in Prop.\ 2.5 also generate nonlinear solutions by superpositions of explicit solutions, namely arbitrary superpositions of eigenmodes with the same Laplacian eigenvalue (the same eigenvalue arises when the wave length is the same), see also \cite{Walsh}. 
In comparison, the superpositions of transverse flows \eqref{eq 3} we present here are much more restricted in the directions, since the wave vectors and flow directions have to satisfy certain orthogonality conditions. 
However, the superposed solutions that we present in this paper exist for any space dimension, and more importantly, these may be superposed modes with different wave lengths \emph{and} (certain) different directions, i.e., with different eigenvalues, and additionally arbitrary wave shapes in case of the Euler equations.\newline
Another example for explicit solutions obtained by superpositions of several flows is presented in \cite{Kambe86}. These consist of a straining flow superposed with several shear flows. 
Two shear layers can merge to a single one or cancel each other out for growing time, depending on the shear layers to be parallel or antiparallel. 
These differ from the solutions and superpositions that we have presented, which on the one hand even do not interact at all in the case of transverse flows \eqref{eq 3}. 
On the other hand, in the case of interacting transverse flows \eqref{eq 6} the plane waves relate to each other by wave length and temporal growth-rate, but without merging or canceling.

The well known parallel flows (see \cite{Wang89} for a general form for the Navier-Stokes equations) consist of a horizontally dependent vertical flow,
which also yields solutions to the rotating Boussinesq equations as shown in \S\ref{s:vanish1}. However, the horizontal plane flows \eqref{sol: vanishBoussinesq} that we identify are differently oriented in the Boussinesq case and depend on a one-dimensional phase variable.

\cite{Weinbaum67} present generalized Beltrami flows to the Navier-Stokes equations, which are two-dimensional solutions with cross flow or axially symmetric solutions with swirl. They also allow non-constant coefficients for their solutions. Our transverse flows \eqref{eq 3} for $n=3$ and $N=2$ belong to these two-dimensional solutions with cross flow.
\cite{Weinbaum67} also present cross flow that solve a heat equation forced by a constant pressure gradient.
We investigate more general forcing, which can also affect the pressure (see Remark \ref{r:forcedGrad}).

For the unforced Navier Stokes equations, analogous to Prop.\ 2.6 in \cite{Majda02}, all these solutions decay exponentially and are thus explicit linear subspaces in the basin of attraction of the trivial solution.
With forcing terms as identified in \S\ref{s:forcing}, nonlinear steady state solutions can be generated via linear variation of constants; a simple case of such forcing is contained, e.g., in \cite{Meshalkin61}.\newline

For the 3D Euler equations, Prop.\ 2.10 in \cite{Majda02} formulates a sufficient Beltrami condition, which our 3D solutions do not satisfy; the Beltrami flow for the Navier-Stokes equations is investigated, for instance, in \cite{Drazin06, Wang89}. As mentioned, for $n\leq 3$ our solutions are \emph{generalized} Beltrami flows \citep{Drazin06,Wang89,Wang90}. We note the other class of  so-called \emph{extended} Beltrami flows of the Navier-Stokes equations \citep{Dyck19}, to which our solutions do not belong.

As mentioned in Remark \ref{r:supParallelFlow}, for the 3D incompressible Navier-Stokes equations \cite{Chai20} derive explicit generalized Beltrami flow solutions with the helical decomposition, which are explicit steady solutions in the 3D Euler equations. These solutions can be understood as a superposition of interacting transverse flows \eqref{eq 6} and adjusted parallel flows. In Remark \ref{r:supParallelFlow} we additionally note the integral superposition and their appearance for $n>3$.

For the Navier-Stokes equations in $\R^n$, $n\geq 2$, and more general linear operator that may also act as forcing, \cite{Chae96} find travelling wave-like solutions with single wave vector. For $n=3$ and usual viscosity \cite{Beloshapkin89} use a superposition of wave modes with different wave vectors but same wave length for steady forcing and solutions. In contrast, the superpositions we allow for can be in integral form and for different wave lengths, but restricted directions. We also allow for different type of forcing as discussed in \S\ref{s:forcing}.

Concerning the inviscid rotating 3D Boussinesq equations, \cite{Majda03} finds solutions with \emph{unbounded} spatially linear velocities and pressure (Theorems 2.4, 2.7)
and also various plane wave type solutions in the non-rotating case. 
With rotation and including the viscous case, the horizontal plane flows \eqref{sol: vanishBoussinesq} we present are barotropic and geostrophically balanced Rossby-type waves, while the interacting horizontal plane flows \eqref{sol: gradBoussinesq} are unbalanced, as are general solutions for the same Laplacian eigenvalue.
Similar, but different from these, are the unbalanced MGW \citep{Yau04, Achatz06}, as discussed in \S\ref{s:vanish1}.

Lastly, we note that explicit generalized Beltrami flows need not be of (superposed) plane wave form, in particular the explicit Lamb-Oseen vortices are the parabolic self-similar solutions \citep{GallayWayne,GohWayne}.\newline

To the best of our knowledge, new results in the present paper are the identification of explicit superposed plane wave type flows with arbitrary wave lengths in (certain) different directions, the explicit representation of solutions with gradient nonlinearity and the corresponding pressure, and the investigation of these explicit solutions with a more general type of adjusted forcing. In an upcoming paper, we will show how solutions of the type discussed in this paper occur in rotating shallow water equations with hyperviscosity, which arises in so-called backscatter modelling, and we discuss stability properties.

Generally, stability properties are a relevant direction of further investigation. We expect that a large part of the steady solutions resulting from those presented are unstable, in particular for large amplitude factors.
For the viscous Boussinesq equations and the pure Navier-Stokes equations, some stability properties of the Kolmogorov flow have been investigated, cf.\ \cite{BalmforthYoung2005} and the references therein. For the inviscid rotating Boussinesq equations various stability studies for the aforementioned  MGW are available in the literature, cf.\ \cite{Achatz06} and the references therein. 
A related line of future investigation concerns the study of modulation equations with respect to the various parameters, similar to what has been done for MGW by a WKB approach, e.g., \cite{Achatz06}.

\bigskip\bigskip
\textbf{Acknowledgements:}
This paper is a contribution to the project M2 (Systematic multi-scale modelling and analysis for geophysical flow) of the Collaborative Research Centre TRR 181 ``Energy Transfers in Atmosphere and Ocean" funded by the Deutsche Forschungsgemeinschaft (DFG, German Research Foundation) under project number 274762653.
The authors thank Marcel Oliver and Ulrich Achatz for fruitful discussions, and Nolan Dyck for comparing with the results in \cite{Dyck19}.\newline

\bibliography{References}
\bibliographystyle{apalike}

\end{document}